\documentclass[12pt]{article}
\usepackage{amsmath, amscd, amssymb,latexsym, epsfig, color}
\input xy
\xyoption{all}

\newtheorem{theorem}{Theorem}[section]
\newtheorem{proposition}[theorem]{Proposition}

\newtheorem{definition}[theorem]{Definition}
\newtheorem{conjecture}[theorem]{Conjecture}

\newtheorem{lemma}[theorem]{Lemma}

\def\proof{\smallskip\noindent {\it Proof: \ }}
\def\endproof{\hfill$\square$\medskip}

\def\R{\mathbb{R}}

\def\M{\mathcal{M}}

\newcommand{\lk}{\mbox{\upshape lk}\,}

\newcommand{\field}{{\bf k}}

\newcommand{\im}{\mbox{\upshape Im}\,}

\newcommand{\Char}{\mbox{\upshape char}\,}

\title{Applications of Klee's Dehn-Sommerville relations.}
\author{Isabella Novik
\thanks{Research partially supported by Alfred P.~Sloan Research
Fellowship and NSF grant DMS-0500748}\\
\small Department of Mathematics, Box 354350\\[-0.8ex]
\small University of Washington, Seattle, WA 98195-4350, USA,\\[-0.8ex]
\small \texttt{novik@math.washington.edu}
\and Ed Swartz
\thanks{Research partially supported by NSF grant DMS-0600502}\\
\small Department of Mathematics, \\[-0.8ex]
\small Cornell University, Ithaca NY, 14853-4201, USA, \\[-0.8ex]
\small \texttt{ebs22@cornell.edu } }

\begin{document}
\maketitle

\begin{abstract}
We use Klee's Dehn-Sommerville relations and other results
on face numbers of homology manifolds without boundary to 
(i) prove Kalai's
conjecture providing lower bounds on the $f$-vectors of
an even-dimensional manifold with all but the middle Betti 
number vanishing, (ii) verify K\"uhnel's conjecture that gives 
an upper bound on the middle Betti number of a $2k$-dimensional 
manifold in terms of $k$ and the number of vertices, and (iii)
partially prove K\"uhnel's conjecture providing upper
bounds on other Betti numbers of odd- and even-dimensional
manifolds. For manifolds with boundary, we derive
an extension of Klee's Dehn-Sommerville relations 
and strengthen Kalai's result on the number of their edges. 

\end{abstract}

\section{Introduction}
In this paper we study face numbers of triangulated manifolds 
(and, more generally, homology manifolds) with and without boundary.
Here we discuss our results deferring most of definitions to subsequent 
sections.

Our starting point is a beautiful theorem known as the Dehn-Sommerville
relations. It asserts that the upper half of the 
face vector of a triangulated manifold 
without boundary  is determined by its Euler 
characteristic together with the lower half of the face vector. 
In this generality the theorem is due to Vic Klee \cite{Klee}. 

Perhaps the most elegant way to present the Dehn-Sommerville relations
is via the $h$-vector of a manifold. The entries of this
vector are certain (alternating)
linear combinations of the face numbers. On the level of $h$-vectors,
the Dehn-Sommerville relations for triangulated 
spheres and odd-dimensional manifolds
merely state that the $h$-vector of these complexes is symmetric.
In the case of spheres, the components of the $h$-vector are 
also known to be positive as they equal dimensions of 
algebraically determined nonzero vector spaces \cite[Chapter 2]{St96}.

Motivated by Dehn-Sommerville relations
together with several commutative algebra results
on  Stanley-Reisner rings of triangulated manifolds, 
Kalai suggested \cite[Section 7]{N98}  
a  modification of the 
$h$-vector, the $h''$-vector, as the ``correct'' $h$-vector
for (orientable) manifolds without boundary. 
The $h''$-vector of orientable manifolds 
(both odd-dimensional and even-dimensional) has since been 
shown to be symmetric \cite{N98}  and nonnegative \cite{NovSw}.

Our first result  is an extension of
Klee's Dehn-Sommerville
relations to manifolds with boundary. Specifically,
we show that for a triangulated manifold with 
a fixed boundary $\Gamma$, the upper half
of the $h$-vector is determined by  
the Euler characteristic, its lower half, and the  
$h$-vector of $\Gamma$. This result is not entirely new. In the
language of  $f$-vectors it was first worked out  
by Macdonald \cite{Macdonald}, and then rediscovered by
Klain \cite{Klain}, and Chen and Yan \cite{ChenYan}. 
However its $h$-vector form appears 
to be absent from the literature. We then use this result
to define a suitable version
of the $h''$-vector for manifolds with boundary 
as well as  show that 
it is symmetric and nonnegative.

Our next result concerns new inequalities 
on the face numbers and Betti numbers
of manifolds without boundary.
Kalai conjectured (private communication) that
the face numbers of a $2k$-dimensional manifold
with all but the middle Betti number vanishing are
simultaneously minimized by the face numbers of a certain 
neighborly $2k$--dimensional manifold. We verify this conjecture.
We also prove a part of a conjecture by 
K\"uhnel \cite[Conjecture 18]{Lutz} that provides an upper
bound on the middle Betti number of a $2k$-dimensional manifold
in terms of $k$ and the number of vertices. Both results turn
out to be a simple consequence of the Dehn-Sommerville relations
and results from \cite{NovSw}.

K\"uhnel further conjectured \cite[Conjecture 18]{Lutz} 
an upper bound on the $i$-th Betti number 
(for all $i$) of a 
$(d-1)$-dimensional manifold with $n$ vertices 
in terms of $i$, $d$, and $n$. We prove that this conjecture
is implied by the $g$-conjecture for spheres. In particular,
K\"uhnel's conjecture holds for manifolds 
all of whose vertex links are polytopal.

In the last section we return to discussing manifolds with boundary.
Here we derive a strengthening of Kalai's theorem 
\cite[Theorem 1.3]{Kalai87} that provides a lower bound on the
number of edges of a manifold in terms of its dimension,
total number of vertices, and the number of interior vertices. Our new
bound also depends on the Betti numbers of the boundary.

The structure of the paper is as follows. In Section 2 we review
necessary background material. In Section 3 we derive the Dehn-Sommerville
relations and define the $h''$-vector for manifolds with boundary.
In Section 4 we deal with Kalai's and K\"uhnel's conjectures.
Finally, in Section 5 we prove a new lower bound
on the number of edges of manifolds with boundary.

\section{Simplicial complexes and face numbers}
In this section we review necessary background material on 
simplicial complexes, Dehn-Sommerville relations, and Stanley-Reisner rings 
of homology manifolds. We refer our readers to  
\cite[Chapter 2]{St96} and the recent paper \cite{NovSw}
for more details on the subject.

Recall that a {\em simplicial complex} $\Delta$ on the vertex set 
$[n]=\{1,2, \ldots, n\}$ is a collection of subsets of $[n]$ that 
is closed under inclusion and contains all singletons $\{i\}$ for 
$i\in[n]$. The elements of $\Delta$ are called {\em faces}. The 
maximal faces (with respect to inclusion) are called {\em facets}.
The {\em dimension of a face} $F\in \Delta$ is $\dim F:=|F|-1$ and
the {\em dimension of $\Delta$} is the maximal dimension of its 
faces. For a simplicial complex $\Delta$ and its face $F$, the 
{\em link} of $F$ in $\Delta$, $\lk(F)$, is the subcomplex of 
$\Delta$ defined by
$$
\lk(F)=\lk_\Delta(F):=\{G\in \Delta \, | \, G\cap F=\emptyset
 \mbox{ and } G\cup F \in \Delta\}.
$$
In particular, the link of the empty face is the complex itself.

A basic combinatorial invariant of a simplicial complex
$\Delta$ on the vertex set $[n]$ is its {\em $f$-vector},
$f(\Delta)= (f_{-1}, f_0, \ldots, f_{d-1})$. Here, $d-1=\dim \Delta$
and $f_i$ denotes the number of $i$-dimensional faces of $\Delta$. Thus
$f_{-1}=1$ (there is only one empty face) and $f_0=n$.
An invariant that contains the same information as the $f$-vector,
but sometimes is more convenient to work with, is the {\em $h$-vector} 
of $\Delta$, $h(\Delta)=(h_0, h_1, \ldots, h_d)$ whose 
entries are defined by the following relation:
\begin{equation}   \label{h-vector}
\sum_{i=0}^d h_i \lambda^{i} = 
\sum_{i=0}^d f_{i-1} \lambda^i(1-\lambda)^{d-i}.
\end{equation}

A central object of this paper is a {\em homology manifold} 
(over a field \field), that is, a $(d-1)$-dimensional pure simplicial 
complex  $\Delta$ such that for all $\emptyset \neq F \in \Delta$,
 the reduced simplicial homology $\tilde{H}_i(\lk F; \field)$ 
vanishes if $i<d-|F|-1$ and is isomorphic to $\field$ or $0$ if  
$i=d-|F|-1$.  A complex is {\em pure} if all of its facets have 
the same dimension.
The {\em boundary faces} of $\Delta$ are those faces $F \neq \emptyset$ 
such that  $\tilde{H}_{d-|F|-1}(\lk F; \field)=0.$ When $\Delta$ has
 no boundary faces,  we write $\partial \Delta = \emptyset$ and $\Delta$ 
is called a homology manifold without boundary.  Otherwise, 
$\partial \Delta$ is the set of boundary faces together with the empty set. 
We will assume that $\partial \Delta$ is a $(d-2)$-dimensional homology 
manifold without boundary.   Under certain conditions this assumption is 
superfluous, see, for instance \cite{Mitchell}.    As demonstrated by the 
suspension of the real projective plane whose `boundary' would be the two 
suspension points for any field whose characteristic is not two, some  
additional assumption is required. We say that $\Delta$ is {\em orientable} 
if the pair $(\Delta, \partial \Delta)$ satisfies the usual 
Poincar\'e-Lefschetz duality associated with orientable compact manifolds 
with boundary.   The prototypical example of a homology manifold 
(with or without boundary) is a triangulation of a topological manifold 
(with or without boundary).

A beautiful theorem due to Klee \cite{Klee} asserts that if 
$\Delta$ is a homology manifold without boundary, then the $f$-numbers 
of $\Delta$ satisfy  linear relations known as the 
{\em Dehn-Sommerville relations}:
\begin{equation} \label{Dehn-Sommerville}
h_{d-i}-h_i = (-1)^{i} {d \choose i}
\left((-1)^{d-1}\tilde{\chi}(\Delta)-1 \right) \quad 
 \mbox{ for all } 0\leq i \leq d.
\end{equation}
Here $\tilde{\chi}(\Delta):=\sum_{i=-1}^{d-1} (-1)^i f_i$ is the 
{\em reduced Euler characteristic} of $\Delta$. Proofs of several results 
in this paper rely heavily on Klee's formula (\ref{Dehn-Sommerville})
and its variations, while other results are concerned with deriving analogs 
of this formula for manifolds with boundary.

In addition to the Dehn-Sommerville relations we exploit several 
results on the Stanley-Reisner rings of homology manifolds.
If $\Delta$ is a simplicial complex on $[n]$, then
its {\em Stanley-Reisner ring} (also called the {\em face ring}) is
$$\field[\Delta] := \field[x_1, \ldots, x_n]/I_\Delta,
     \quad \mbox{where }
I_\Delta=(x_{i_1}x_{i_2}\cdots x_{i_k} : \{i_1<i_2<\cdots<i_k\}
\notin\Delta).
$$
(Here and throughout the paper $\field$ is an infinite field of an
arbitrary characteristic.)
Since $I_\Delta$ is a monomial ideal, the ring $\field[\Delta]$ is 
graded, and we denote by $\field[\Delta]_i$ its $i$th homogeneous 
component. The Hilbert series of $\field[\Delta]$, 
$F(\field[\Delta], \lambda):= 
\sum_{i=0}^\infty \dim_{\field} \field[\Delta]_i \cdot \lambda^i$, 
has the following properties.

\begin{theorem} \label{Stanley} {\rm {\bf (Stanley)}} 
Let $\Delta$ be a $(d-1)$-dimensional simplicial complex. Then
$$
F(\field[\Delta], \lambda) = \frac{\sum_{i=0}^d h_i \lambda^i}{(1-\lambda)^d}.
$$
\end{theorem}

\begin{theorem} \label{Schenzel}
{\rm {\bf (Schenzel)}} Let $\Delta$ be a $(d-1)$-dimensional 
homology manifold, and let $\theta_1, \ldots, \theta_d \in \field[\Delta]_1$ 
be such that 
$\field[\Delta]/\Theta:=\field[\Delta]/(\theta_1, \cdots, \theta_d)$ 
is a finite-dimensional vector space over \field. Then 
$$
F(\field[\Delta]/\Theta, \lambda) = 
\sum_{i=0}^d \left(h_i(\Delta)+
  {d \choose i} \sum_{j=1}^{i-1} (-1)^{i-j-1} \beta_{j-1}(\Delta)\right) 
\cdot \lambda^i,
$$
where $\beta_{j-1}:=\dim_\field \tilde{H}_{j-1}(\Delta; \field)$.
\end{theorem}

Theorem \ref{Stanley} can be found in \cite[Theorem II.1.4]{St96},
while Theorem \ref{Schenzel} is from \cite{Schenzel}.  In view of 
Theorem \ref{Schenzel}, for a $(d-1)$-dimensional
homology manifold $\Delta$, define 
\begin{equation} \label{h'}
h'_i(\Delta):=h_i(\Delta)+
{d \choose i} \sum_{j=1}^{i-1} (-1)^{i-j-1} \beta_{j-1}(\Delta).
\end{equation}
We remark that if $|\field|=\infty$ then a set of linear forms 
$\{\theta_1, \ldots, \theta_d\}$ satisfying the assumptions of 
Theorem \ref{Schenzel} always exists, e.g., choosing 
``generic'' $\theta_1, \ldots, \theta_d$ does the job .

The following theorem summarizes several results 
on the $h'$-numbers of 
homology manifolds that will be needed later on. 
For $0<m={x \choose i}:=x(x-1)\cdots(x-i+1)/i!$ where
$0<x\in\R$, define $m^{<i>}:={x+1 \choose i+1}$. Also set $0^{<i>}:=0$.

\begin{theorem} \label{h'-properties}
Let $\Delta$ be a $(d-1)$-dimensional homology manifold. Then
\begin{enumerate}
\item $h'_0=1$, $h'_1=f_0-d$, and for all $1\leq i \leq d$
$$
h'_i\geq {d \choose i}\beta_{i-1} \qquad \mbox{and} \qquad
h'_{i+1} \leq \left(h'_i-{d \choose i}\beta_{i-1} \right)^{<i>}.
$$
\item Moreover, if $\Delta$ is a homology manifold without 
boundary that is orientable over \field, i.e., 
$\beta_{d-1}(\Delta)=\beta_0(\Delta)+1$, then
\begin{equation} \label{h'-DS}
h'_{d-i}-h'_i ={d \choose i} (\beta_i-\beta_{i-1}) \qquad \mbox{ for all } 
0\leq i\leq d.
\end{equation}
\end{enumerate}
\end{theorem}
Part 1 of this theorem was recently proved in 
\cite{NovSw} (see Theorems 3.5 and 4.3 there).
Part 2 is a simple variation of Klee's 
Dehn-Sommerville relations, see \cite[Lemma 5.1]{N98}. 
It is obtained by combining equations 
(\ref{Dehn-Sommerville}) and  (\ref{h'})
with Poincar\'e duality for homology manifolds.

Eq.~(\ref{Dehn-Sommerville}) 
implies that all homology spheres and odd-dimensional
manifolds without boundary satisfy $h_i=h_{d-i}$ for all $i$. 
While this symmetry fails for even-dimensional manifolds 
with $\tilde{\chi}\neq 1$, Theorem \ref{h'-properties} together 
with Poincar\'e duality suggests we consider the following modification
of the $h$-vector and yields the following algebraic version of 
(\ref{Dehn-Sommerville}).

\begin{proposition} \label{h''}
Let $\Delta$ be a $(d-1)$-dimensional 
homology manifold without boundary. 
Assume further that $\Delta$ is connected and orientable over \field.
Let 
$$h''_d:=h'_d \mbox{ and } 
h''_i(\Delta):= h'_i(\Delta) -{d \choose i}\beta_{i-1}(\Delta) =
h_i-{d \choose i}\sum_{j=1}^{i} (-1)^{i-j}\beta_{j-1}, 
\mbox{ for } 0\leq i \leq d-1.$$
Then $h''_i\geq 0$  and
$h''_i(\Delta)=h''_{d-i}(\Delta)$ for all $0\leq i \leq d$.
\end{proposition}

In view of Proposition \ref{h''} and results of \cite{NovSw2}
that interpret $h''$-numbers as dimensions of homogeneous components 
of a  Gorenstein ring, $h''$ can be regarded as the ``correct" 
$h$-vector for orientable homology manifolds without boundary. 
What is the analog of $h''$ for manifolds with boundary? 
We deal with this question in the following section.

\section{Dehn-Sommerville for manifolds with boundary}
Klee's equations (\ref{Dehn-Sommerville}) generate
a complete set of linear 
relations satisfied by the $h$-vectors of 
homology manifolds with empty boundary. More generally,
one can fix a (non-empty) homology manifold 
$\Gamma$ and ask for the 
set of all linear relations satisfied by the $h$-vectors of 
homology manifolds whose boundary is $\Gamma$. 
Deriving such relations and defining what seems to be the 
``correct'' version of the $h''$-vector is the goal of this section.

We have the following version of Dehn-Sommerville relations: 
\begin{theorem}  \label{DS-boundary}
Let $\Delta$ be a $(d-1)$-dimensional homology manifold with
boundary. Then 
$$
h_{d-i}(\Delta)-h_i(\Delta)=
{d \choose i}
(-1)^{d-1-i}\tilde{\chi}(\Delta) -
g_i(\partial\Delta) \quad \mbox{ for all } 0\leq i \leq d,
$$
where 
$g_i(\partial\Delta):=h_i(\partial\Delta)-h_{i-1}(\partial\Delta)$.
\end{theorem}
\proof Write $f_i:=f_i(\Delta)$ and $h_i:=h_i(\Delta)$. 
Let $f^b_i:=f_i(\partial\Delta)$, and define 
$h^b_i$ and $g^b_i$ in a similar way. 
Also let $f^\circ_i:=f_i(\Delta)-f^b_i$ be the 
``interior" $f$-vector, and let $h^\circ_i$ be defined from $f^\circ$ 
according to Eq.~(\ref{h-vector}). With this notation, 
we obtain from \cite[Corollary II.7.2]{St96} that
$$
(-1)^d F(\field[\Delta], 1/\lambda) = 
(-1)^{d-1}\tilde{\chi}(\Delta) + 
\sum_{i=1}^{d}\frac{f^\circ_{i-1} \lambda^i}{(1-\lambda)^i}.
$$
Substituting Theorem \ref{Stanley} in
the above formula yields
$$
(-1)^d\frac{\sum_{i=0}^d h_{d-i}\lambda^i}{(\lambda-1)^d}=
(-1)^{d-1}\tilde{\chi}(\Delta) +
\sum_{i=1}^{d}
  \frac{f^\circ_{i-1} \lambda^i (1-\lambda)^{d-i}}{(1-\lambda)^d}
 =(-1)^{d-1}\tilde{\chi}(\Delta) +
\frac{\sum_{i=0}^{d} h^\circ_i \lambda^i}{(1-\lambda)^d}, $$
which is  equivalent to
$$
\sum_{i=0}^d (h_{d-i}-h^\circ_i)\lambda^i=
(-1)^{d-1}\tilde{\chi}(\Delta)(1-\lambda)^d.$$
Subtracting $\sum_{i=0}^d g^b_i \lambda^i$ from both sides and 
noting that $h^\circ_i+g^b_i=h_i$, implies the result.
\endproof

\smallskip\noindent While Theorem \ref{DS-boundary} 
appears to be new, $f$-vector forms of the same equality have appeared before.
   Chen and Yan gave a generalization which applies to more general stratified 
spaces \cite{ChenYan}.  However, we believe that the first place where an 
equivalent formula appears is due to Macdonald \cite{Macdonald}.

We now turn to finding the right definition of 
$h''$ for orientable homology manifolds with boundary.
Recall that a connected $(d-1)$-dimensional 
homology manifold $\Delta$
is orientable over \field \ if 
$H_{d-1}(\Delta, \partial\Delta; \field)\cong \field$. 
By Poincar\'e-Lefschetz duality, 
if $\Delta$ is such a manifold, then
$H_{i-1}(\Delta, \partial\Delta)\cong H_{d-i}(\Delta)$. 
Write $\beta_{i-1}(\Delta, \partial\Delta)$
to denote
$\dim H_{i-1}(\Delta, \partial\Delta)$.

We start by expressing $g_i(\partial\Delta)$ in terms of its 
Betti and $h'$-numbers. Substituting Eq.~(\ref{h'}) 
in $g_i(\partial\Delta) =
h_i(\partial\Delta)-h_{i-1}(\partial\Delta)$ and recalling that
$\dim(\partial\Delta)=d-2$, we obtain
\begin{equation}  \label{g}
g_i(\partial\Delta) = \left[ h'_i(\partial\Delta)- h'_{i-1}(\partial\Delta)
+{d-1 \choose i-1} \beta_{i-2}(\partial\Delta)\right]
+{d \choose i} \sum_{j=1}^{i-1}
 (-1)^{i-j}\beta_{j-1}(\partial\Delta).
\end{equation}

\begin{theorem} \label{h''_boundary}
Let $\Delta$ be a $(d-1)$-dimensional homology manifold with 
nonempty boundary. If $\Delta$ is orientable, then for all $0\leq i < d$,
$$
h'_{d-i}(\Delta) -{d \choose d-i}\beta_{d-i-1}(\Delta)=
h'_i(\Delta)-\overline{g}_i(\partial\Delta)-
{d \choose i} \dim \im(H_{i-1}(\Delta) \stackrel{\psi}{\rightarrow} 
H_{i-1}(\Delta, \partial\Delta)),
$$
where $\overline{g}_i(\partial\Delta):=
h'_i(\partial\Delta)- h'_{i-1}(\partial\Delta)
+{d-1 \choose i-1} \beta_{i-2}(\partial\Delta)$ and $\psi$ 
is the map in the long exact sequence of the pair $(\Delta, \partial \Delta).$
\end{theorem}

\proof If $i=0$, then both sides are equal to 0. 
For $0< i <d$, 
using Eq.~(\ref{h'}) and Theorem \ref{DS-boundary},
we obtain
\begin{eqnarray*}
h'_{d-i}(\Delta) &-& {d \choose d-i}\beta_{d-i-1}(\Delta) \\
&=& 
h_i(\Delta)-g_i(\partial\Delta) 
+ (-1)^{d-1-i}{d \choose i}
\left[\tilde{\chi}(\Delta)
+\sum_{j=1}^{d-i}(-1)^{j}\beta_{j-1}(\Delta)\right] \\
&=& h_i(\Delta)-g_i(\partial\Delta) 
+ (-1)^{d-1-i}{d \choose i}
\left[\sum_{j=d-i+1}^{d} (-1)^{j-1}\beta_{j-1}(\Delta)\right]\\
&=& h_i(\Delta)-g_i(\partial\Delta) 
+ (-1)^{i}{d \choose i} 
\sum_{j=0}^{i-1}(-1)^{j}\beta_{j}(\Delta, \partial\Delta),
\end{eqnarray*}
where the last step is by Poincar\'e-Lefschetz duality.
Substituting equations (\ref{h'}) and (\ref{g}) in the
last expression then yields,
\begin{eqnarray*}
h'_{d-i}(\Delta) &-& {d \choose d-i}\beta_{d-i-1}(\Delta) \\
&=&
h'_i(\Delta)-\overline{g}_i(\partial\Delta) 
-{d \choose i} \sum_{j=0}^{i-1}(-1)^{j-i-1} 
\left[\beta_{j}(\Delta, \partial\Delta)
-\beta_{j-1}(\partial\Delta)+\beta_{j-1}(\Delta)\right].
\end{eqnarray*}
The result follows, since by long exact homology sequence of 
the pair $(\Delta, \partial\Delta)$, the last summand equals 
$-{d \choose i}\dim \im(H_{i-1}(\Delta) \rightarrow 
H_{i-1}(\Delta, \partial\Delta))$.
\endproof 

Theorem \ref{h''_boundary} 
suggests the following definition of the $h''$-vector and 
shows (together with theorem \ref{h'-properties}) that it is symmetric
and non-negative.
\begin{definition} For $\Delta$ --- a $(d-1)$-dimensional orientable 
homology manifold with a nonempty boundary, define
$$
h''_i(\Delta)=\left\{
\begin{array}{ll} h'_i(\Delta)-\overline{g}_i(\partial\Delta)-
{d \choose i} \dim \im(H_{i-1}(\Delta) \rightarrow 
H_{i-1}(\Delta, \partial\Delta))& \quad \mbox{ for } 
 i\leq d/2 \\
h'_i(\Delta)-{d \choose i}\beta_{i-1}(\Delta)
 & \quad \mbox{ for } 
 i > d/2.
\end{array}
 \right.
$$
\end{definition}
Note that in the case of the empty boundary and $i<d$,
this definition agrees with the one given in 
Proposition \ref{h''}.

\section{Manifolds without boundary: Kalai's and K\"uhnel's conjectures}
In this section we settle a conjecture of Kalai that 
provides lower bounds for the face numbers of even-dimensional 
homology manifolds with all Betti numbers but the middle one vanishing. We
also partially settle a conjecture by K\"uhnel on the Betti numbers 
of homology manifolds. Throughout this section, $\Delta$ 
denotes a $(d-1)$-dimensional orientable homology manifold without boundary. 
Note that if $\field$ is a field of characteristic two, then this class 
includes all triangulated topological manifolds without boundary.

We start by discussing even-dimensional manifolds. 
The following result was conjectured by K\"uhnel 
\cite[Conjecture 18]{Lutz}.

\begin{theorem}  \label{Kuhnel-middle}
Let $\Delta$ be a $2k$-dimensional orientable homology manifold 
with $n$ vertices. Then 
$${2k+1 \choose k}\beta_k(\Delta) \leq {n-k-2 \choose k+1}.
$$
Moreover, if equality is attained then $\beta_i=0$ for all $i<k$.
\end{theorem}
\proof
Choose a nonnegative real number $x$ such that 
$$
h'_{k}-{2k+1 \choose k}\beta_{k-1} ={x \choose k}.
$$ 
It exists since according to Theorem \ref{h'-properties}, 
$h'_{k}-{2k+1 \choose k}\beta_{k-1}\geq 0$. 
Moreover, the same theorem implies that
$h'_{k+1}\leq {x+1 \choose k+1}$. Thus
$$
{2k+1 \choose k}\beta_k \stackrel{\mbox{\small{by (\ref{h'-DS})}}}{=}
 h'_{k+1}-h'_{k}+
{2k+1 \choose k}\beta_{k-1}\leq {x+1 \choose k+1}-{x \choose k}=
{x \choose k+1}.
$$
Finally, since $h'_1=n-2k-1$, 
another application of Theorem \ref{h'-properties} shows that
$h'_k \leq {n-k-2 \choose k}$, hence $x\leq n-k-2$, 
and 
${2k+1 \choose k}\beta_k \leq {n-k-2 \choose k+1}$, as required.
Furthermore, equality implies that 
$h'_i={n-k-2 \choose i} =(h'_{i-1})^{<i>}$ for all $2 \leq i\leq k+1$, 
which by Theorem \ref{h'-properties} is possible only if $\beta_i=0$ 
for all $i<k$.
\endproof
 
Theorem \ref{Kuhnel-middle} implies that if $\beta_k\geq 1$,
then $n-k-2 \geq 2k+1$, or equivalently, $n\geq 3k+3$. In other words, 
having a non-vanishing middle Betti number requires at
least $3k+3$ vertices. (This result was originally proved by  Brehm and 
K\"uhnel for PL-triangulations \cite{BrKuh}.) Moreover, if such a 
homology manifold,  $\M_k$, has exactly $3k+3$ vertices, then 
$$
\beta_k (\M_k)=1, \, \beta_i(\M_k)=0 \mbox{ for }i<k, \, 
\mbox{ and } h_i(\M_k)=h'_i(\M_k)={k+1+i \choose i} 
\mbox{ for }i\leq k+1.
$$
In particular, the face numbers of $\M_k$ (whether it exists or not) 
are uniquely determined by Eqs.~(\ref{h-vector}) 
and (\ref{Dehn-Sommerville}).  These face numbers turn out to be minimal 
in the following sense (as was conjectured by 
Gil Kalai, personal communication):

\begin{theorem} \label{Kalai's-conj}
Let $\Delta$ be a $2k$-dimensional orientable homology manifold
with $\beta_k\neq 0$ being the only non vanishing Betti number out of
all $\beta_l$, $l\leq k$. Then 
$$f_{i-1}(\Delta)\geq f_{i-1}(\M_k) \quad \mbox{for all}
\quad  1\leq i \leq 2k+1.
$$
\end{theorem}

\proof Substituting  $\beta_l=0$, $l<k$, in 
Theorem \ref{Schenzel} and Eq.~(\ref{Dehn-Sommerville}), we obtain that 
$$
h_j(\Delta)= h'_j(\Delta)  \mbox { and }
h_{k+j+1}(\Delta)= h_{k-j}(\Delta) + 
(-1)^j {2k+1 \choose k-j}\beta_k(\Delta) \quad \mbox{for } 0\leq j\leq k.
$$
Eq.~(\ref{h-vector}) then implies 
\begin{eqnarray}
f_{i-1}(\Delta) &=&
\sum_{j=0}^i {2k+1-j \choose 2k+1-i} h'_j(\Delta), \quad
      \mbox{ if } i\leq k,  \quad \mbox{ and } \label{small_i}\\
f_{i-1}(\Delta) &=& 
\sum_{j=0}^k 
\left[ {2k+1-j \choose 2k+1-i} + {j \choose 2k+1-i} \right]
 h'_j(\Delta) \nonumber \\
& + & \beta_k(\Delta)\left[\sum_{j=0}^{i-k-1} 
 (-1)^{j}{k-j \choose 2k+1-i}{2k+1 \choose k-j}\right] \quad
 \mbox{ if $i\geq k+1$}.  \label{big_i}
\end{eqnarray}
Since (i) the same formulas apply to the $f$-numbers of $\M_k$, 
(ii) the coefficients of the $h'$-numbers in 
Eqs.~(\ref{small_i}) and (\ref{big_i}) are nonnegative, 
and (iii) $\beta_k(\Delta)\geq 1 =\beta_k(\M_k)$, to complete the proof 
it only remains to show that $h'_i(\Delta)\geq h'_i(\M_k)$
for all $i\leq k$ and that the coefficient of $\beta_k$ in 
Eq.~(\ref{big_i}) is nonnegative for all $i\geq k+1$. 

The latter 
assertion follows by noting that the sequence 
$$a_j= {k-j \choose 2k+1-i}{2k+1 \choose k-j}, \quad 0\leq j \leq i-k-1
$$
is decreasing (indeed, 
$a_j/a_{j+1} = (k+2+j)/(i-k-1-j)>1$), and hence 
$a_0-a_1+\cdots +(-1)^{i-k-1}a_{i-k-1} \geq 0$. 

To verify the former assertion, we use the same trick as in the proof of 
Theorem \ref{Kuhnel-middle}. Let $0\leq x \in \R$ 
be such that $h'_k(\Delta)={x \choose k}$. Then according to 
Theorem \ref{h'-properties}, 
$h'_{k+1}(\Delta)\leq {x+1 \choose k+1}$ while
$h'_{k+1}(\Delta)-h'_k(\Delta) ={2k+1 \choose k+1}\beta_k 
\geq {2k+1 \choose k+1}$. Thus we have
$$
{2k+1 \choose k+1} \leq h'_{k+1}(\Delta)-h'_k(\Delta) \leq 
{x+1 \choose k+1} -{x \choose k}={x \choose k+1}.
$$
Hence $x\geq 2k+1$, and so $h'_k(\Delta)\geq {2k+1 \choose k}$.
Applying Theorem \ref{h'-properties} once again, we infer that
$h'_i(\Delta)\geq {k+1+i \choose i}=h'_i(\M_k)$ 
for all $i\leq k$.
\endproof

In addition to Theorem \ref{Kuhnel-middle}, K\"uhnel conjectured 
(see \cite[Conjecture 18]{Lutz}) that a $(d-1)$-dimensional
manifold with $n$ vertices satisfies
${d+1 \choose j+1}\beta_j(\Delta)\leq {n-d+j-1 \choose j+1}$ 
for all  $0\leq j \leq \lfloor d/2 \rfloor -1$. The case of 
$j=0$ merely says that every connected component of $\Delta$ 
has at least $d+1$ vertices. The case of $j=1$ is equivalent to 
Kalai's lower bound conjecture \cite[Conjecture 14.1]{Kalai87}
that was recently settled in \cite[Theorem 5.2]{NovSw}. 
For other values
of $j$ we have the following partial result.
We recall that a $(d-1)$-dimensional homology sphere $\Gamma$ 
is said to have the {\em hard Lefschetz property} if for 
a generic choice of 
$\theta_1, \ldots, \theta_d, \omega \in \field[\Gamma]_1$,
the map 
$$\field[\Gamma]/(\theta_1, \ldots, \theta_d)_i 
\stackrel{\cdot\omega^{d-2i}}{\longrightarrow} 
\field[\Gamma]/(\theta_1, \ldots, \theta_d)_{d-i}$$
is an isomorphism of $\field$-spaces for all $i\leq d/2$.
It is a result of Stanley \cite{St80} that in the case of 
$\Char \field=0$ all simplicial 
polytopes have this property, and it is the celebrated 
$g$-conjecture that all homology spheres do.
 
\begin{theorem} \label{g->Kuhnel}
Let $\Delta$ be a $(d-1)$-dimensional 
orientable homology manifold with $n$ vertices. If for
every vertex $v$ of $\Delta$ the link of $v$ has the 
hard Lefschetz property (e.g., $\Char \field=0$ and all vertex links
are polytopal spheres), then 
$${d+1 \choose j+1}\beta_j(\Delta)\leq {n-d+j-1 \choose j+1}
\quad \mbox{
for all } 0\leq j \leq \lfloor \frac{d}{2} \rfloor -1.
$$ 
If equality is attained for some $j=j_0$, then $\beta_i=0$ for all
$i\neq j_0$, $0\leq i \leq \lfloor d/2 \rfloor -1$.
\end{theorem}

\proof 
Since all vertex links of $\Delta$ 
have the hard Lefschetz property, 
Theorem 4.26 of \cite{Sw} implies that for 
a sufficiently generic choice of 
$\theta_1,\ldots, \theta_d, \omega \in \field[\Delta]_1$ and every 
$j\leq \lfloor d/2 \rfloor -1$, the linear map 
$$
\field[\Delta]/(\theta_1, \ldots, \theta_d)_{d-j-1} 
\stackrel{\cdot\omega}{\longrightarrow} 
\field[\Delta]/(\theta_1, \ldots, \theta_d)_{d-j}$$
is surjective. The dimensions of the spaces involved
are $h'_{d-j-1}$ and $h'_{d-j}$, respectively 
(see Theorem \ref{Schenzel}). Also, by \cite[Cor.~3.6]{NovSw},
the dimension of the kernel of this map is at least 
${d \choose d-j-1} \beta_{d-j-2}$. Therefore,
\begin{equation}  \label{surj}
h'_{d-j} \leq h'_{d-j-1}-{d \choose d-j-1} \beta_{d-j-2}
\quad \mbox{ for all }j\leq \lfloor d/2 \rfloor -1.
\end{equation}
Apply Poincar\'e duality and Eq.~(\ref{h'-DS}) to rewrite this
inequality in the form
$$
h'_j+{d \choose j}(\beta_j - \beta_{j-1}) 
\leq h'_{j+1}- {d \choose j+1}\beta_j,
$$ or, equivalently,
\begin{equation} \label{h'<=h'}
{d+1 \choose j+1}\beta_j 
\leq h'_{j+1}-\left[h'_j - {d \choose j}\beta_{j-1}\right].
\end{equation}

Let $0\leq x\in \R$ be such that 
$h'_{j+1}={x+1 \choose j+1}$. 
Then by Theorem \ref{h'-properties}, 
$h'_j - {d \choose j}\beta_{j-1} \geq {x \choose j}$, 
and so the right-hand-side of (\ref{h'<=h'}) is 
$\leq {x \choose j+1}$. Also, since $h'_1=n-d$, 
$h'_{j+1}\leq {n-d+j \choose j+1}$, and hence 
$x\leq n-d+j-1$. Thus
${d+1 \choose j+1}\beta_j \leq {x \choose j+1} 
\leq {n-d+j-1 \choose j+1}$, as required.

If equality occurs for some
$j=j_0$, then $x=n-d+j_0-1$, and we obtain that
$h'_{i+1}={n-d+i \choose i+1}=(h'_i)^{<i>}$ for 
all $i\leq j_0$. By Theorem \ref{h'-properties} 
this can happen only if $\beta_{i-1}=0$ for all $i\leq j_0$. 
Moreover in this case, $h_{i+1}={n-d+i \choose i+1}$ 
for all $i\leq j_0$, hence $\Delta$ is $(j_0+1)$-neighborly
(that is, every set of $j_0+1$ vertices of $\Delta$ is a face of
$\Delta$).

What about $\beta_i$ for $i>j_0$? To prove that all 
these Betti numbers vanish as well, note that for 
equality ${d+1 \choose j_0+1}\beta_{j_0}
= {n-d+j_0-1 \choose j_0+1}$ to happen, the
inequality in (\ref{surj}) should hold as equality for $j=j_0$.
The same argument as in the proof of \cite[Theorem 5.2]{NovSw}
then shows that $h_{j_0}(\lk v)=h_{j_0+1}(\lk v)$ for every 
vertex $v$ of $\Delta$. Since, by our assumptions, all vertex
links of $\Delta$ satisfy the $g$-conjecture, and since
$\Delta$ is $(j_0+1)$-neighborly, we conclude that for every 
vertex $v$,
 $h_i(\lk v)=h_{j_0}(\lk v)$ for all 
$j_0\leq i \leq (d-1)/2$, and that $h(\lk v)=h(\lk w)$ for 
all vertices $v$ and $w$ of $\Delta$.
This information about links turns out to be
 enough to compute the entire $h$-vector of $\Delta$. Indeed, it follows 
from \cite[Remark 4.3]{HerNov} that
$$
h_r(\Delta)=(-1)^r{d \choose r} + 
\sum_{i=0}^{r-1} (-1)^{r-i-1}
 \frac{(d-1-i)!i!}{(d-r)!r!}\cdot n \cdot h_i(\lk v).
$$
Hence 
$$
g_{r+1}(\Delta)=h_{r+1}-h_r = {d+1 \choose r+1}
\left[ (-1)^{r+1} + 
n \sum_{i=0}^{r-1}  
   \frac{(-1)^{r-i} h_i (\lk v)}{(d-i){d \choose i}} 
+ \frac{n\cdot h_r(\lk v)}{(r+1){d+1 \choose r+1}} \right],
$$
and since $h_{j_0}(\lk v)= h_{j_0+1}(\lk v)=\cdots$, 
we infer
that for all $j_0+1\leq r \leq (d-1)/2$,
$$
\frac{g_{r+1}}{{d+1 \choose r+1}} +
\frac{g_r}{{d+1 \choose r}} 
= n\cdot h_{j_0}(\lk v) 
\left[ -\frac{1}{(d-r+1){d \choose r-1}} +
\frac{1}{(r+1){d+1 \choose r+1}} + \frac{1}{r{d+1 \choose r}}
\right]= 0.
$$
Therefore, 
$$
(-1)^{r-j_0}\frac{g_{r+1}}{{d+1 \choose r+1}}=
\frac{g_{j_0+1}}{{d+1 \choose j_0+1}}=
\frac{h_{j_0+1}-h_{j_0}}{{d+1 \choose j_0+1}} =\beta_{j_0}.
$$
Substituting this result in 
Eq.~(\ref{h'<=h'}) with $j=j_0+1$  and using that all
$\beta_i$ for $i<j_0$ vanish, yields
${d+1 \choose j_0+1} \beta_{j_0+1}\leq 
\left[h_{j_0+2}+{d \choose j_0+2}\beta_{j_0}\right]
-\left[h_{j_0+1}-{d \choose j_0+1}\beta_{j_0}\right]
=g_{j_0+2}+{d+1 \choose j_0+2}\beta_{j_0} =0,
$
and so $\beta_{j_0+1}=0$. Assuming by induction that 
$\beta_{j_0+1}=\ldots=\beta_{r-1}=0$, a similar computation 
using Eq.~(\ref{h'<=h'}) with $j=r$ then 
implies that $\beta_r=0$ for all $j_0<r\leq (d-1)/2$.
\endproof

Kalai conjectured \cite[Conj.~14.2]{Kalai87} that if $\Delta$ 
is a $(d-1)$-dimensional manifold without boundary, then
$
h''_{j+1}(\Delta)-h''_j(\Delta) 
\geq {d \choose j}\beta_j(\Delta).
$
This is an immediate consequence of
Eq.~(\ref{h'<=h'}). Thus Kalai's conjecture holds 
for all manifolds whose vertex links have the hard Lefschetz 
property.

\section{Rigidity inequality for manifolds with boundary}
In this section we return to our discussion of the face numbers of 
homology manifolds with nonempty boundary. The goal here is to strengthen 
Kalai's result \cite[Theorem 1.3]{Kalai87} asserting that if $\Delta$ is a 
$(d-1)$-dimensional manifold with boundary and $d\geq 3$, then 
$h_2(\Delta)\geq f_0^\circ (\Delta)$, where as in Section 3, 
$f_0^\circ (\Delta)$ denotes the number of interior vertices of $\Delta$.
Our main result is

\begin{theorem}  \label{h2-boundary}
If $\Delta$ is a connected $(d-1)$-dimensional homology manifold with 
nonempty orientable boundary and $d\geq 5$, then 
\begin{equation} \label{h2-eq}
h_2(\Delta)\geq f_0^\circ (\Delta)+{d \choose 2}\beta_1(\partial\Delta) + 
d\ \beta_0(\partial \Delta).
\end{equation}
If $d=4$ and the characteristic of \field \ is two, then
$$h_2(\Delta)\geq f_0^\circ (\Delta)+3\  \beta_1(\partial\Delta) + 
4 \ \beta_0(\partial \Delta).$$
\end{theorem}

Since the boundary of a $3$-manifold with boundary is a collection of closed 
surfaces, using a field whose characteristic is two maximizes the relevant 
Betti numbers, so we have restricted ourselves to this case.  
Before beginning the proof of this theorem we establish 
some preliminary results pertaining to rigidity in characteristic 
$p > 0.$  In characteristic zero the cone lemma, gluing lemma, 
and Proposition \ref{h_2+g_2}  follow easily from the work of
 Kalai \cite{Kalai87} and Lee \cite{Lee}.

\begin{definition}
A $(d-1)$-dimensional complex $\Delta$ is $\field$-{\bf rigid} 
if for generic $\theta_1, \dots, \theta_{d+1}$ linear forms
and $1 \le i \le d+1,$ multiplication 
$\cdot \theta_i: \field[\Delta]/(\theta_1, \dots, \theta_{i-1})_1 \to  
\field[\Delta]/(\theta_1, \dots, \theta_{i-1})_2$ is injective. 
\end{definition}
It follows from Proposition \ref{h_2+g_2} below that if $\Delta$ is
\field-rigid, then the
$\field$-dimension of $\field[\Delta]/(\theta_1, \dots, \theta_{d})$ 
is $h_2(\Delta)$ and of $\field[\Delta]/(\theta_1, \dots, \theta_{d+1})$ 
is $g_2(\Delta)$. In fact, it is not hard to see (but we will not use it
here) that the converse holds as well.

\begin{lemma} (Cone lemma)
If $\Delta$ is \field-rigid, then the cone on 
$\Delta, C(\Delta),$ is \field-rigid.
\end{lemma}

\proof  Observe that 
$\field[C(\Delta)]\cong\field[\Delta]\otimes_\field \field[x_0]$. Hence 
for any $\theta_0$ of the form $x_0+\sum_{i=1}^n \alpha_ix_i$, 
$\theta_0$ is a non-zero-divisor on $\field[C(\Delta)]_1$ and
the quotient ring
$\field[C(\Delta)]/(\theta_0)$ is isomorphic to $\field[\Delta]$. The
assertion follows.
\endproof

\begin{lemma} (Gluing lemma)
If $\Delta_1$ and $\Delta_2$ are $(d-1)$-dimensional \field-rigid complexes 
and there are at least $d$ vertices in $\Delta_1 \cap \Delta_2,$ then 
$\Delta_1 \cup \Delta_2$ is \field-rigid.
\end{lemma}

\proof
Set $\Delta = \Delta_1 \cup \Delta_2.$  Since $\Delta_l$ ($l=1,2$) 
is a subcomplex of $\Delta$, there is a natural surjection 
$\field[\Delta] \longrightarrow \field[\Delta]_l$.
Consider the following commutative square.
\begin{equation}  \label{commsq}
\begin{CD} 
\field[\Delta]/(\theta_1, \dots, \theta_{i-1})_2 @>>> 
\field[\Delta_1]/(\overline{\theta}_1, 
\dots,\overline{\theta}_{i-1})_2 @>>> 0\\
@AA\cdot \theta_iA @AA\cdot \theta_iA \\
\field[\Delta]/(\theta_1, \dots,\theta_{i-1})_1 @>>> 
\field[\Delta_1]/(\overline{\theta}_1, \dots,\overline{\theta}_{i-1})_1 @>>> 0.
\end{CD}
\end{equation}
Here, $\overline{\theta}$ is the image of $\theta$ in $\field[\Delta_1].$ 
Suppose $\omega$ is in the kernel of the left-hand vertical map.  
Then its image in $\field[\Delta_1]/(\theta_1, \dots, \theta_{i-1})_1$ 
must be in the kernel of the right-hand vertical map, and hence zero when 
restricted to $\field[\Delta_1]/(\theta_1, \dots, \theta_{i-1})_1.$  
Similarly, $\omega$ is zero when restricted to 
$\field[\Delta_2]/(\theta_1, \dots, \theta_{i-1})_1.$  
But, if there are at least $d$ vertices in $\Delta_1 \cap \Delta_2$,
 then $\omega = 0$ in $\field[\Delta]/(\theta_1, \dots, \theta_{i-1})_1.$
\endproof

\begin{proposition} \label{h_2+g_2}
Let $\Delta_1, \dots, \Delta_b$ be \field-rigid $(d-1)$-dimensional complexes 
with disjoint sets of vertices.  If $\Delta = \cup \Delta_i,$ then for 
generic linear forms $\Theta = (\theta_1, \dots, \theta_d)$ and 
$\omega$, 
$$\dim_\field \left(\field[\Delta]/\Theta\right)_2 = 
h_2(\Delta) + {d \choose 2} (b-1)
\, \mbox{ and }
\, \dim_\field \ker
\left[(\field[\Delta]/\Theta)_1 \stackrel{\cdot\omega}{\to} 
(\field[\Delta]/\Theta)_2\right]
=d (b-1).$$
\end{proposition}

\proof
Suppose that $w$ is in  the kernel of 
$\cdot \theta_1: \field[\Delta]_1 \to \field[\Delta]_2$.  
Using a commutative square analogous to (\ref{commsq}), 
we see that restricted to each vertex set $w$ is zero.  
Hence $w=0.$  Therefore,
 $$\dim_\field \left(\field[\Delta] /(\theta_1)\right)_2 =
 \dim_\field \field[\Delta]_2 - \dim_\field \field[\Delta]_1 = 
(f_1 + f_0) - f_0 = f_1.$$
 
 Now replace $\field[\Delta]$ with $\field[\Delta]/(\theta_1)$ 
in (\ref{commsq}) and consider multiplication by $\theta_2.$  
The same argument shows that any $w$ in the kernel must restrict 
to a multiple of $\theta_1$ on the vertex set of each $\Delta_j.$  
The dimension of the space of such $w$ in 
$(\field[\Delta]/(\theta_1))_1$ is $b-1.$  Thus, 
$$\dim_\field (\field[\Delta]/(\theta_1, \theta_2))_2 = 
f_1 - (f_0-1) + (b-1).$$  
Continuing with this reasoning we see that for each $i$ 
the dimension of the kernel of multiplication by $\theta_i$ 
on $(\field[\Delta]/(\theta_1, \dots, \theta_{i-1}))_1$ is $(i-1) (b-1).$  
Hence, for $i \ge 2,$
 $$\dim_\field (\field[\Delta]/(\theta_1, \dots, \theta_{i-1}))_2 = 
f_1 - (i-2) f_0 + {i-1 \choose 2} + {i-1 \choose 2} (b-1).$$
Setting $i = d+1$ finishes the proof.
\endproof

\smallskip\noindent {\it Proof of Theorem \ref{h2-boundary}: \ }  
First we consider the situation when $d \ge 5.$
Let   $\Gamma$ be the simplicial complex obtained from 
$\Delta$ by coning off each component of the boundary of $\Delta.$  
Specifically, let $c_1, \dots, c_b$ be the components of the boundary of 
$\partial \Delta.$  We introduce new vertices 
${\bf n+1}, \dots, {\bf n+b}$ and set
$$
\Sigma=(({\bf n+1}) \ast c_1) \cup \dots \cup (({\bf n+b}) \ast c_b)
\quad \mbox{ and } \quad
\Gamma = \Delta \cup \Sigma.$$
Then $\Gamma$ is a $(d-1)$-dimensional pseudomanifold  
that is \field-rigid.   The proof is by induction on $d.$  
Any $\Delta$ homeomorphic to $S^2$ is $\field$-rigid. This follows from 
\cite[Cor.~3.5]{Murai}. 
So the cone lemma implies that the closed star of a vertex in 
a three-dimensional \field-homology sphere is 
\field-rigid.  Now, using the gluing lemma we can take the union with  
closed stars of other vertices until we see that an arbitrary
 three-dimensional \field-homology sphere is \field-rigid.  
Since for every vertex $v\in\Delta$, the link of $v$ in $\Gamma$ is a 
\field-homology
sphere, induction on $d$ implies that this link is $\field$-rigid. Hence
the closed star of $v$ in $\Gamma$ is $\field$-rigid for all $v\in \Delta$. 
Taking the union of the closed stars of the noncone points 
using the gluing lemma shows that $\Gamma$ is \field-rigid.

Observe that $f_0(\Gamma)=f_0(\Delta)+b$ and 
$f_1(\Gamma)=f_1(\Delta)+f_0(\partial \Delta)$. Thus 
$h_2(\Gamma)=h_2(\Delta)+h_1(\partial\Delta) - (d-1)\beta_0(\partial \Delta)$. 
For $\Sigma$ we have $f_0(\Sigma) = f_0(\partial \Delta) + b$ and 
$f_1(\Sigma) = f_1(\partial \Delta) + f_0(\partial \Delta).$  
Hence, $h_2(\Sigma) = h_2(\partial \Delta)  - (d-1) \beta_0(\partial \Delta).$

Consider the face rings $\field[\Gamma]$ and $\field[\Sigma]$, 
and let 
$\theta_1, \ldots, \theta_d,\omega\in \field[\Gamma]_1$ be generic 
linear forms. Since $\Sigma$ is a subcomplex of $\Gamma$, 
there is a natural
surjection 
$\phi: \field[\Gamma] \longrightarrow \field[\Sigma]$. 
Let $\overline{\theta}_i$  denote the image of $\theta_i$  
under $\phi$, and consider 
$\field(\Gamma):=\field[\Gamma]/(\theta_1, \ldots, \theta_d)$ and
 $\field(\Sigma):=
\field[\Sigma]/(\overline{\theta}_1, \ldots, 
\overline{\theta}_d)$. Then $\phi$ induces a surjection
$\field(\Gamma) \longrightarrow \field(\Sigma)$.
Denoting by $I \subset \field(\Gamma)$ its kernel, we obtain the following 
commutative diagram whose rows are exact:
\begin{equation} \label{commd}\minCDarrowwidth22pt\begin{CD}
0 @>>>   I_2
 @>>> \field(\Gamma)_2
 @>>> \field(\Sigma)_2 @>>> 0 \\
@. @AA{\cdot\omega}A  @AA{\cdot\omega}A  @AA{\cdot\omega}A\\
0 @>>>   I_1 
 @>>> \field(\Gamma)_1
 @>>> \field(\Sigma)_1 @>>> 0. \\
\end{CD}
\end{equation}

Since $\Gamma$ is \field-rigid, $\dim \field(\Gamma)_2=h_2(\Gamma)$ and 
the middle vertical map is an injection.
Hence the left vertical map is also an injection.  
By the cone lemma and the argument which proved that 
$\Gamma$ is \field-rigid,  each of the $b$ components of $\Sigma$ is 
\field-rigid.  Proposition \ref{h_2+g_2} says that 
$\dim_\field \field(\Sigma)_2 = 
h_2(\Sigma) + {d \choose 2} \beta_0(\partial \Delta).$

By  Proposition \ref{h_2+g_2}, the dimension of the kernel of the right 
vertical map  is  $d \beta_0(\partial \Delta).$   Applying the snake lemma, 
we find that the dimension of the cokernel of $\cdot \omega: I_1 \to I_2$ 
is at least $d \beta_0(\partial \Delta)$ and thus 
$\dim I_1 + d \beta_0(\partial \Delta)  \le \dim I_2.$

What are the dimensions of $I_1$ and $I_2$?
From exactness of rows, we infer that
\begin{equation}  \label{dimI1}
\dim I_1= \dim \field(\Gamma)_1-\dim \field(\Sigma)_1=
(f_0(\Delta)+b-d)-(f_0(\partial\Delta)+b-d)=f_0^\circ(\Delta),
\end{equation}
and
\begin{eqnarray}
\dim I_2&=& \dim \field(\Gamma)_2-\dim \field(\Sigma)_2=
 \nonumber\\
&=& h_2(\Gamma)-h_2(\Sigma)-{d \choose 2}\beta_0(\partial\Delta) \nonumber \\
& =&h_2(\Delta)- g_2(\partial\Delta) - {d \choose 2}\beta_0(\partial\Delta)
\nonumber \\
&\leq& h_2(\Delta)-
\left[{d \choose 2}\beta_1(\partial\Delta)-
{d \choose 2}\beta_0(\partial\Delta)\right]-
{d \choose 2}\beta_0(\partial\Delta) \nonumber\\
&=& h_2(\Delta)-{d \choose 2}\beta_1(\partial\Delta),  \label{dimI2}
\end{eqnarray}
where the penultimate step follows from \cite[Theorem 5.2]{NovSw} applied to
connected components of $\partial\Delta$ and from the observation that
for a $(d-2)$-dimensional 
complex $\partial\Delta$,  its $g_2$-number equals the sum of the 
$g_2$-numbers of its connected components minus 
${d \choose 2}\beta_0(\partial\Delta)$. Comparing the right-hand-sides of
(\ref{dimI1}) and (\ref{dimI2}) and using $\dim I_1 + d\beta_0 \leq \dim I_2,$
implies the result.

Two modifications are necessary when $d=4.$   First, each component of the 
boundary of $\Delta$ is a closed surface,  so the Dehn-Sommerville relations 
tell us that the $g_2$ of each component is $3 \beta_1.$   Second,   to show 
that $\Sigma$ is \field-rigid the induction must begin with any closed 
surface instead of just $S^2.$  In his thesis \cite{Fogelsanger}, 
Fogelsanger proved that any triangulation of a closed surface is generically 
$3$-rigid in the graph-theoretic sense.  Fogelsanger used three properties of 
generic 3-rigidity:  a cone lemma, a gluing lemma, and a result of Whiteley's 
concerning vertex splitting \cite{Whiteley}.  Our cone lemma and gluing lemma 
cover the first two.   Whiteley's vertex splitting result, combined 
with \cite[Theorem 10]{Lee} due to Carl Lee,  is characteristic independent. 
   Hence, Fogelsanger's proof shows that a triangulation of a closed surface 
is \field-rigid. \endproof

  Now we give a series of examples that show that for any 
$d \ge 5, \beta_1, \beta_0$ and $f^\circ_0$,  
Theorem~\ref{h2-boundary} is optimal.  We recall a family  of complexes 
introduced by K\"uhnel and Lassman.
  
  \begin{theorem} \cite{KuhLas}
  For every $d \ge 4$ and $n \ge 2d-1$ there exists a complex 
$M^d(n)$ with $n$ vertices such that
  \begin{itemize}
   \item
     $M^d(n)$ is a  $B^{d-2}$-bundle over the circle.    
In particular, $M^d(n)$ is a manifold with boundary.
     
     \item Depending on the parity of $n$ and $d$  
the boundary of $M^d(n)$ is either $S^{d-2} \times S^1$ 
 or the nonorientable $S^{d-2}$-bundle over the circle.  
Hence, for $d \ge 5,$ the  first Betti number of 
$\partial M^d(n)$ is one for  any field.  
When $d =4$ and the characteristic of \field \  is $2,$ 
then $\beta_1(\partial M^4(n)) = 2.$
     
     \item
     $h_2(M^d(n)) = {d \choose 2}.$
     
     \item
     All of the vertices are on the boundary of $M^d(n).$  
The link of every vertex is combinatorially equivalent to a stacked polytope.  
  
  \end{itemize}
  \end{theorem}
  
  Evidently, $M^d(n)$ for $d \ge 5$ is an example of equality in 
Theorem \ref{h2-boundary} with $\beta_1(\partial \Delta) = 1$ 
and $f^\circ_0 = 0.$  For spaces with $\beta_1(\partial \Delta) > 1,$ 
begin with two disjoint copies of $M^d(n).$ 
Choose two $(d-2)$-faces on their respective boundaries and a bijection 
between their vertices.  Now identify these vertices and associated faces 
according to the chosen bijection.  The resulting space has no interior
 vertices and is a manifold with boundary whose boundary is topologically 
the connected sum of two copies of the boundary of $M^d(n).$  
Thus the first Betti number is now two.  Direct computation shows that
 $h_2$ of the new space is $2 {d \choose 2}.$  Repeating this operation of 
connected sum along the boundary $b$ times  with $M^d(n)$ produces an example 
of equality in Theorem \ref{h2-boundary} with $\beta_1(\partial \Delta) = b$ 
 and $f^\circ_0 = 0.$  To construct $\Delta$ with $f^\circ_0 = m > 0$ 
simply take a complex with $f^\circ_0=0$ and subdivide a facet $m$ times.  
Each such subdivision increases $h_2$ and $f^\circ_0$ by one while leaving 
the topological type of the complex unchanged.    
  
  To produce spaces $\Delta$ with $\beta_0(\partial 
  \Delta) > 0$, begin with any of the above examples.   
It is possible to subdivide a facet $d$ times so that there 
is now a facet with interior vertices.  See \cite{CSS} for the algorithm. 
 Removing the open facet leaves a manifold whose  boundary  has two components,
 the original and the boundary of the simplex.  The new space will have the 
same number of interior vertices and its $h_2$ will have increased by $d.$  

In dimension three, the same constructions lead to examples of equality in 
Theorem~\ref{h2-boundary} with  arbitrary $f^\circ, \beta_0,$ and even 
$\beta_1.$  Since the boundary of a three-dimensional manifold $\Delta$ must 
have even Euler characteristic, this is the best we can hope for.

 All of the complexes constructed using the procedures have the property 
that the link of every boundary vertex is combinatorially a stacked polytope, 
and the link of every interior vertex is a stacked sphere. 
  
 \begin{conjecture}
 If $\Delta$ is a connected $(d-1)$-dimensional homology manifold with 
nonempty orientable boundary and $d\geq 4$, then 
equality occurs in Theorem \ref{h2-boundary} if and only if all of the links of $\Delta$ are combinatorially equivalent to 
stacked polytopes or stacked spheres.
\end{conjecture}

\end{document}